\input amstex
\documentstyle{amsppt}
\magnification = \magstep1
\hsize=6.2truein
\vsize=8.4truein

\define\al{\alpha}
\define\be{\beta}
\define\ga{\gamma}

\define\de{\delta}

\define\ep{\epsilon}
\define\la{\lambda}
\define\om{\omega}
\define\si{\sigma}
\define\ta{\tau}

\define\ct{{\Cal T}}
\define\cf{{\Cal F}}
\topmatter
\title
Patterns of Dependence among Powers of Polynomials
\endtitle
\author
Bruce Reznick
\endauthor
\date
June 3, 2001
\enddate
\address
Department of Mathematics,
University of Illinois at Urbana-Champaign,
1409 W. Green Street,
Urbana, IL 61801
\endaddress
\email
reznick\@math.uiuc.edu
\endemail
\abstract
Let $F = \{f_1,\dots,f_r\}$ be a family of polynomials and let the
ticket of $F$, $T(F)$,
denote the set of integers $m$ so that $\{f_j^m\}$ is linearly
dependent. We show that $|T(F)| \le \binom {r-1}2$ and present many
concrete examples.
\endabstract
\endtopmatter

\document

\head
1. Introduction and definitions
\endhead

The motivation for this  paper is a remarkable
example due to A. H. Desboves in 1880 (see [3,\ p.684])  and,
independently, N. Elkies [2,\ p.542] in 1995.   
\proclaim{Example 1}
\endproclaim
Let
$$
\gathered
 f_1(x,y) = x^2 + \sqrt{2} xy - y^2, \quad  f_2(x,y) =i x^2 - \sqrt{2}
 xy +i y^2 , \\ f_3(x,y) =  -x^2 +
\sqrt{2} xy + y^2, \quad  f_4(x,y) = -ix^2 -
\sqrt{2} xy -i y^2.
\endgathered
\tag 1.1
$$
Then $\sum_j f_j^5 = 0$. It is not hard to show that $\sum_j
f_j =\sum_j f_j^2 = 0$ as well, but it can be shown that $\{f_j^m\}$
is linearly independent for other integers $m$. This
is surprising, and not only because of the gap in the exponents. Let
$M_m(t_1,t_2,t_3,t_4)= \sum_j t_j^m$. Then 
 $(f_1,f_2,f_3,f_4)$ parameterizes the intersection
$\{M_1 = 0\} \cap \{M_2 = 0\}\cap \{M_5 = 0\} \subset \bold C^4$.
This is not what one would expect geometrically; however, 
Newton's Theorem for symmetric forms implies that $M_5 \in
\left(M_1,M_2\right)$, so $ \{M_5 = 0\} \subset \bigl(\{M_1 = 0\} \cap \{M_2
= 0\}\bigr)$.    We shall return to (1.1) throughout
this paper, and derive it in several different ways.  

\smallskip
Let $F = \{f_j\}$ be a finite set of  polynomials.
 The {\it ticket} of $F$, $T(F)$, is defined by 
$$
T(F) = \{m \in \bold N: \{f_j^m\} \text{ is linearly
dependent} \}. 
$$
We let $\overline{m}(F)$ denote the maximum element of $T(F)$, with the
understanding that $\overline{m}(F) = 0$ if $T(F) = \emptyset$. More
specifically, for integers $r \ge 2$, $n \ge 2$,
$d\ge 1$, let $\cf(r,n,d)$ denote the 
set of families $F = \{f_1,\dots, f_r\}$, where each $f_j$ is a homogeneous
polynomial over $\bold C$ in $n$ variables of degree $d$, and no two
$f_j$'s are proportional. 
Let 
$$
\ct(r,n,d) = \{T(F): F \in \cf(r,n,d)\}.
$$
Let $\overline{m}_r$ denote the maximum of $\overline{m}(F)$, given
$|F| = r$. It is known that $\overline{m}_r \le r^2-2r$.
Finally, in view of the geometry of Example 1, we  say that  $F \in
\cf(r,n,d)$ is a {\it dysfunctional} family if $|T(F)| > r-2$, the term is
used in the sense of excessive dependence.
\proclaim
{Example 2}
\endproclaim
Here are some families with tickets that are easy to verify by hand.
$$
\gather
T(\{x,y,z\}) = \emptyset, \tag 1.2)(i \\
T(\{x,y,x+y\})  = \{1\}, \tag 1.2)(ii \\
T(\{x^2-y^2,2xy,x^2+y^2\}) = \{2\}. \tag 1.2)(iii
\endgather
$$
\smallskip
One of the principal goals of this research is to determine
$\ct(r,n,d)$. We are able to do so in the cases $r \le 3$, $(n,d) =
(2,1)$,  and $(r,n,d) = (4,2,2)$. Besides specifying $(r,n,d)$, we can
consider the sets which appear as the ticket of any family. Let
$$
\ct := \bigcup_{r,n,d} \ct(r,n.d).
$$
The following  tends to strike
listeners as either very  plausible or ridiculous.

\proclaim{Conjecture  1}
$\ct$ consists of all finite subsets of $\bold N$.
\endproclaim

There does not seem to be much in the literature on this subject,
except for $\overline{m}_r$. If $r=3$, then after scaling, $m \in
T(F)$ implies that $f_1^m + f_2^m - f_3^m = 0$. In 1879,
J. Liouville proved (see [13,\ p.263]) that Fermat's Last Theorem is true for
non-constant polynomials, hence $\overline{m}_3 = 2$.

M. Green proved in 1975 [4,\ p.71] that if  $\{\phi_j\}$, $1
\le j \le r$, are holomorphic functions in $n$ complex variables, no two
of which are proportional, and $\sum_{j=1}^r \phi_j^m = 0$, then $m
\le (r-1)^2-1$. This implies that  $\overline{m}_r \le r(r-2)$, and hence
 $|T(F)| \le r(r-2)$.

D\. J\. Newman and M\. Slater proved in 1979  [10,\ p.481] that a non-zero
sum of the $m$-th powers of $r$ non-constant polynomials in a
single complex variable has degree  $\ge m - \binom
r2$. (Their proof goes through for $n$ variables as well.) By taking
this sum to be 1 and homogenizing, it follows that if $F \in \cf(r,n,d)$
and one of the $f_j$'s is a $d$-th power of a linear form, then
$\overline{m}(F) \le \binom r2 - 1$. They also show [10,\
p.486] that $\overline{m}_r \le 8(r-1)^2 -1$ for 
unrestricted $F$, but this is weaker than Green's bound. A 1972 example
due to J. Molluzzo [9] shows that  $\overline{m}_r \ge \lfloor 
\frac {(r-1)^2}4\rfloor$. Molluzzo's theorem is discussed in Example
12 below; the bound is improved here for $r \le 8$. 

There is a connection to Nevanlinna theory. If $m \in T(F)$
for $F \in \cf(r,2,d)$, then by dehomogenizing and transposing the
linear relation among $f_j^m$, we obtain an
equation $\sum_{j=1}^{r-1}\phi_j^m(z) = 1$ in rational functions
$\phi_j$.  A recent survey of
this work has been made by W. Hayman [7]; see also [8]. 
G. Gundersen  [5] has found
meromorphic (not rational) functions $g_j$ for which $g_1^6(z) + g_2^6(z)+
g_3^6(z) = 1$.

\smallskip
We collect some trivial remarks, which will be used without citation
throughout the rest of the paper. 
\smallskip
  If  $f_j = \al f_i$, then $f_i^m$
and $f_j^m$ are dependent for 
every $m$,  so $T(F) = \bold N$ is an uninteresting ticket. This 
explains the non-proportionality condition in the definition of
$\cf(r,n,d)$. 

The dimension of the vector space of forms of degree $D$ in $n$
variables is $\binom {n+D-1}{n-1}$. Thus, if $F \in \cf(r,n,d)$ and
if $r > \binom {n+md-1}{n-1}$, then $\{f_j^m\}$ must be dependent. In
this case we say that $m$ is {\it forced} to be in  $T(F)$. 

 There is no {\it a priori} reason to assume that the polynomials
in a given family $F$  are all homogeneous and of the same
degree. Suppose $F = \{f_j\}$, $1 \le j \le r$, where each $f_j$ is a
polynomial in $n$ variables of degree $\le d$. Homogenize $F$ by
defining 
$$
f'_j(x_1,\dots,x_{n+1}) = x_{n+1}^df_j\left(\frac{x_1}{x_{n+1}},\dots,
\frac{x_n}{x_{n+1}}\right),\qquad 1 \le j \le r ,
$$
and let $F' = \{f_j'\} \in \cf(r,n+1,d)$. It is easy to see that  $T(F) =
T(F')$. Similarly, if we
dehomogenize $G = \{g_j\} \in \cf(r,n,d)$, via
$$
g_j'(x_1,\dots,x_{n-1}) = g_j(x_1,\dots,x_{n-1},1),\qquad  1 \le j \le r ,
$$
and let $G' = \{g_j'\}$, then $T(G)=T(G')$. Thus, there is no
loss of generality in restricting ourselves to the various $\cf(r,n,d)$'s.

 If  $F = \{f_j\}$ and
$$
f'_j(x_1,\dots,x_n) = f_j\left(\tsize\sum\limits_{k=1}^n
\al_{1k}x_k + \be_1, 
\dots\tsize\sum\limits_{k=1}^n \al_{nk}x_k+\be_n\right),\qquad   1 \le j \le r ,
$$
where $\det([\al_{jk}]) \neq 0$, and  $F' = \{f_j'\}$, then clearly $T(F) =
T(F')$: the ticket of a family is not changed by a simultaneous
invertible linear change of variables.

 If $F = \{f_j\}$, $f_j' = \la_jf_j$, where $\la_j \neq 0$  and
$F' = \{f_j'\}$, then $T(F) = T(F')$. That is,
we may normalize the members of a family without affecting its ticket.
In particular, if $\sum_j
c_jf_j^m = 0$, $c_j \neq 0$, then we may
assume without loss of generality that $c_j = \pm 1$ at our
convenience. (Of course, we can only do this for {\it one} such
exponent $m \in T(F)$.)

 Given $\{f_j\}$ in $n$ variables and $n$ polynomials
$\phi_k(y_1,\dots,y_m)$, we can define
$$
f_j'(y_1,\dots,y_m) := f_j(\phi_1(y),\dots,\phi_n(y)).
$$
Then $\sum_{j=1}^r \la_j f_j = 0$ implies $\sum_{j=1}^r \la_j f_j' = 0$.
It is possible that  the $f_j'$'s might not be pairwise
non-proportional, but it might also be possible that a judicious
choice of $\phi_k$ will allow for the ticket to be extended.
For example, $T(\{x+y,x-y\}) = \emptyset$ and $4xy = (x+y)^2 - (x-y)^2$ is
not a square, but if we replace $x$ and $y$ by $x^2$ and $y^2$ and add
$2xy$ to $F$, we put 2 into the ticket; c.f., (1.2)(iii).

If $0 \neq g$ is a form, then $T(\{f_j\}) = T(\{gf_j\})$. Thus, $d < d'$
implies that $\ct(r,n,d) \subseteq \ct(r,n,d')$.

If $n < n'$, then by viewing  $f_j(x_1,\dots,x_n)$
as a form in  $(x_1,\dots,x_n,\dots,x_{n'})$, 
we see that  $\ct(r,n,d) \subseteq \ct(r,n',d)$. 

 Regrettably, it is not true that $r < r'$ implies $\ct(r,n,d)
\subseteq \ct(r',n,d)$.
 For example, $T(\{x,y\}) = \emptyset$, so
$\emptyset \in \ct(2,2,1)$. However, $m=1$ is forced to be in $T(F)$
for every $F \in \cf(3,2,1)$. 
(We suspect that the only exceptions to the inclusion $\ct(r,n,d)
\subseteq \ct(r',n,d)$ are consequences of such forcing.)

\smallskip

Here is an outline of the rest of this paper. In section two, we prove
our main theoretical result: if $F \in \cf(r,n,d)$, then $|T(F)| \le
\binom{r-1}2$. The proof is effective, in that  we give an explicit 
polynomial in $(m,x_1,\dots,x_n)$ which vanishes identically when
$m \in T(F)$. This modestly improves Green's $|T(F)| \le r(r-2)$.
 In section
three, we analyze $\ct(r,n,d)$ in the ``easy" cases $r \le 3$, 
and $(n,d) = (2,1)$ and $(r,n,d) = (4,2,2)$, and discuss $(n,d) = (\ge
3,1)$. We also show that $\{1,2,3\} \notin \ct(4,n,d)$. In section
four we work out some concrete tickets. The most dysfunctional family
we have found is
$F \in \cf(2k,2,2)$  for which $|T(F)| = 3k-3$ (and $\overline{m}(F) =
4k-4$.)  For each integer $a$, we give $F_a \in \cf(a+2,2,a)$ so that $T(F_a)$
consists precisely of the divisors of $a$. We also give explicit
$F, F' \in \cf(6,2,2)$ with $T(F) = \{1,2,8\}$ and $T(F') =
\{1,2,3,4,8,14\}$. These generalize (1.1). In fact, for every
$v \ge 2$, there is a non-trivial 
identity equating two sums of the $(3v-1)$-st powers of $v$ binary
quadratic forms. We close the paper with
generalizations, speculations, open questions and acknowledgments.

\head
2. How long can a ticket be?
\endhead

In this section we prove that $|T(F)| \le \binom {r-1}2$. 
The exposition is simplified by assuming that the $f_j$'s are not
homogeneous. We assume that $f_j = f_j(x)$, where $x =
(x_1,\dots,x_n)$, although $n$ does not play an explicit role in the
proof. 
The proof is constructive in that, for each $F$, we define a polynomial
of degree $\binom {r-1}2$ in $(m,x)$ which vanishes identically in $x$
whenever $m \in T(F)$.

\proclaim{Theorem 1}
Suppose $F = \{f_1,\dots,f_r\}$ is a set of polynomials over $\bold C$, no
two of which are proportional, with $\deg(f_j) \le d$. Then $|T(F)|
\le \binom{r-1}2$. 
If $d \le r-3$, then  $|T(F)| < \binom{r-1}2$.
\endproclaim
\demo{Proof}
Suppose first that $ \{g_1,\dots, g_r\}$ are linearly
dependent polynomials and $g_j = \sum_{k \ge 0} 
(g_j)_k$, where each $(g_j)_k$ is homogeneous of degree $k$. Then
$\sum_{j=1}^r \la_jg_j = 0$ implies that 
$\sum_{j=1}^r \la_j(g_j)_k = 0$ for all $k$. That is, 
$$
\left( \matrix (g_1)_0 & (g_2)_0 & \cdots & (g_r)_0 \\
 (g_1)_1 & (g_2)_1 & \cdots & (g_r)_1 \\
\cdots & \cdots & \ddots & \cdots \\
 (g_1)_{r-1} & (g_2)_{r-1} & \cdots & (g_r)_{r-1} \endmatrix \right)
\left( \matrix \la_1 \\ \la_2 \\ \vdots \\ \la_r  \endmatrix\right) =
\left( \matrix 0 \\ 0 \\ \vdots \\ 0  \endmatrix\right).
\tag 2.1
$$
Therefore, as a polynomial in $x$,
$$
W(g_1,\dots,g_r;x): = 
 \vmatrix (g_1)_0 & (g_2)_0 & \cdots & (g_r)_0 \\
 (g_1)_1 & (g_2)_1 & \cdots & (g_r)_1 \\
\cdots & \cdots & \ddots & \cdots \\
 (g_1)_{r-1} & (g_2)_{r-1} & \cdots & (g_r)_{r-1} \endvmatrix = 0.
\tag 2.2
$$

(This Wronskian argument is assuredly not a {\it sufficient} condition
for dependence.  
Given (2.2), the only null-vectors satisfying (2.1) might well have
non-constant components. It might be  the case that deg($g_j) < r- 1$
for all $j$, rendering (2.2) trivial, and it might happen 
that  $\sum_{j=1}^r \la_jg_j$ has terms of degree $\ge r$.
 Nevertheless, the weak necessary condition (2.2) will suffice for
our needs.)

We turn to $F$.  Since $f_i/f_j$ is not 
constant for $i \neq 
j$, there exists $k = k(i,j)$ so that
$$
\frac{\partial}{\partial x_k}\left(\frac{f_i}{f_j}\right) = f_j^{-2} \left(f_j
\frac{\partial f_i}{\partial x_k} - f_i
\frac{\partial f_j}{\partial x_k}\right) \tag 2.3
$$
is not identically zero. For precision, let $k(i,j)$ be the smallest
such index. 
Further, there exists a point $P \in \bold C^n$ which is {\it not}
contained in the finite  union
$$
\left(\bigcup_{j=1}^r \{f_j = 0\}\right)\bigcup\left(\bigcup_{1 \le i
< j \le r}  \left\{f_j
\frac{\partial f_i}{\partial x_{k(i,j)}} - f_i
\frac{\partial f_j}{\partial x_{k(i,j)}} = 0\right\}\right). \tag 2.4
$$ 
By translation, we may assume that $P=0$. If we now let $\ga_j =
(f_j(0))^{-1}$, and replace $f_j$ by $\ga_jf_j$, then we affect neither
the definition of the varieties in (2.4) nor $T(F)$. Accordingly, we may
assume without loss of generality that 
$f_j(0) = 1$ for all $j$. By evaluating (2.3) at $P=0$, we see that
 $\frac{\partial f_i}{\partial
x_{k(i,j)}}(0) \neq \frac{\partial f_j}{\partial x_{k(i,j)}}(0)$, 
hence the first-order terms in the Taylor series 
of the $f_j$'s at 0 are different. To emphasize this point, we write
$(f_j)_1 = L_j$.

We now compute $(f_j^m)_k$ explicitly:
$$
\gathered
f_j^m = \left( 1 + L_j + (f_j)_2 + \dots + (f_j)_d\right)^m \\=
\sum\Sb \ell_0 + \dots + \ell_d = m \\ \ell_0 \ge 0, \dots, \ell_d \ge
0 \endSb \frac {m!}{\ell_0!\dots \ell_d!} 
1^{\ell_0}L_j^{\ell_1}((f_j)_2)^{\ell_2} \dots ((f_j)_d)^{\ell_d},
\endgathered
\tag 2.5
$$
and so $(f_j^m)_k$ will consist of those terms from (2.5) in which 
$\sum_{i=1}^d i\ell_i = k$.  After noting that $\ell_0 = m -
\sum_{i=1}^d \ell_i$ and writing
$(N)_s :=\frac{N!}{(N-s)!} = N(N-1)\cdots(N-(s-1))$, we have 
$$
(f_j^m)_k =
\sum \Sb  
\ell_1 + \dots +\ell_d \le  m \\ \ell_1 + \dots + d\ell_d 
= k \endSb \frac{(m)_{\ell_1 + \dots + \ell_d}}{\ell_1!\dots \ell_d!} 
L_j^{\ell_1}((f_j)_2)^{\ell_2} \dots ((f_j)_d)^{\ell_d}.
\tag 2.6
$$

Consider (2.6) as presenting $(f_j^m)_k$ as a polynomial in $m$,
whose coefficients are polynomials in $x$.
From this perspective, each summand has $m$-degree $\ell_1 + \dots +
\ell_d$. This degree is $\le k$, and the maximum is achieved precisely
when $\ell_1 = k$ and $\ell_i = 0$, $i \ge 2$. The unique term in
$(f_j^m)_k$ with $m$-degree equal to $k$ is 
$\frac{L_j^k}{k!}(m)_k$. Thus,
$$
(f_j^m)_k = \tfrac {L_j^k}{k!}\cdot m^k + {\Cal O}_x(m^{k-1}), \qquad 0 \le
j \le r-1.\tag 2.7
$$
It is also easy to see that the term in $(f_j^m)_k$ with smallest
$m$-degree occurs when $\ell_d = \lfloor k/d \rfloor$, and all other
$\ell_i$'s vanish, except for $\ell_{k - d \lfloor \frac kd \rfloor} =
1$, (unless  $ d \ | \ k$).  Let $u(k,d):= \lfloor \frac{k-1}d
\rfloor$. This smallest degree is $1 + u(k,d)$ and implies that
$$
m(m-1)\cdots(m - u(k,d))\ |\ (f_j^m)_k,\qquad 1 \le k \le r-1
$$
(More directly, 
since $f_j^m$ has degree $\le dm$, $(f_j^m)_k = 0$ if $k > md$; that
is, if $m \le \frac{k-1}d$.) Note that, in particular, $m\ |\
(f_j^m)_k$ if $k \ge 1$. 

We apply (2.2) with $g_j = f_j^m$. Let $W(m;x) :=
W(f_1^m,\dots,f_r^m;x)$. If $m \in T(F)$, then $W(m;x) = 0$ as a
polynomial in $x$. In view of (2.7), we have
$$
W(m;x) = 
\gathered
\vmatrix 1  &  \cdots & 1 \\
  L_1\cdot m &  \cdots &  L_r\cdot  m\\
\frac {L_1^2}2\cdot  m^2 +  {\Cal O}_x(m)  & \cdots & \frac
{L_r^2}2\cdot  m^2 + {\Cal O}_x(m)\\  
\cdots  & \ddots & \cdots \\
 \frac {L_1^{r-1}}{(r-1)! }\cdot m^{r-1} +  {\Cal O}_x(m^{r-2}) & 
  \cdots &
 \frac {L_{r-1}^{r}}{(r-1)!}\cdot  m^{r-1} +  {\Cal O}_x(m^{r-2})
 \endvmatrix. 
\endgathered
\tag 2.8
$$
Then $W(m;x)$ has $m$-degree $\le \binom r2$, and the
coefficient of $m^{\binom r2}$  is the Vandermonde determinant
$$
\tfrac 1{2!\cdots(r-1)!} \vmatrix  1   & \cdots & 1 \\
L_1(x)  & \cdots &  L_r(x) \\
\cdot  & \ddots & \cdots \\
L_1^{r-1}(x)  & \cdots &  L_r^{r-1}(x) \endvmatrix
= \tfrac 1{2!\cdots(r-1)!} \prod_{1 \le i < j \le r}\left( L_j(x) - L_i(x)
 \right). 
$$
Since the $L_j$'s are distinct linear forms, there exists 
$y \in \bold C^n$ so that  $L_i(y) \neq L_j(y)$ for $i \neq j$. Thus,
$W(m;y)$ is a polynomial in $m$ of exact degree $\binom r2$. Since
$m^{r-1}$ is a factor, $W(m;y)$ can have at most $\binom r2 - (r-1) =
\binom{r-1}2$ distinct positive integer roots, and these contain
$T(F)$. 

Moreover, the $k$-th row of (2.8) is divisible by $(m)_{1+u(k-1,d)}$, hence
$W(m;y)$ has a factor of  
$$
\prod_{k=0}^{r-1} (m)_{1+u(k,d)} =
\prod_{k=0}^{r-1} m(m-1)\dots (m-u(k,d)) = 
\prod_{i=0}^{u(r-1,d)}
(m-i)^{r-id-1}.
$$
In particular, if $d \le r - 3$, then $m^{r-1}(m-1)^2$ divides
$W(m;y)$, so $|T(F)| < \binom{r-1}2$. Indeed,  if $F \in
\cf(r,n,d)$, then
$$
|T(F)| \le \binom r2 - (r-1) - \sum_{i=1}^{u(r-1,d)} (r-2-id).
$$
If $d = 2$, this bound works out to be $\lfloor r^2/4\rfloor - 1$.
For larger $d$ it is $\approx \bigl(\frac{d-1}{2d}\bigr) r^2$.
\qed
\enddemo
If $r=3$, then $|T(F)| \le 1$; see Theorem 3 below.

If $r=4$, then  $|T(F)| \le 3$, and the bound is achieved in (1.1).
Example 5 will discuss $F \in \ct(4,2,2)$ with  $T(F) = \{1,2,4\}$.
Theorem 5 will show that there is no family in any $\cf(4,n,d)$ for which
$T(F) = \{1,2,3\}$. (However, by Theorem 2, this is the ticket for any
$F \in \cf(5,2,1)$.)

Theorem 1 gives an algorithm for computing the ticket of a given
family: once a finite set of candidate exponents is found,
the dependence of the $f_j^m$'s can be checked. 

If $(n,d) = (1,2)$ and $r=4$, then (2.5) becomes explicitly 
$$
\gathered
(1 + at + bt^2)^m = 
1 + mat + \left(\tfrac {m(m-1)}2a^2 + mb\right)t^2 \\ + \left( \tfrac
{m(m-1)(m-2)}6 a^3 + m(m-1)ab \right) t^3 + \dots,
\endgathered
\tag 2.9
$$
(Note the factor of $(m)_0$ in $(f^m)_1$ and  $(f^m)_2$, and the
factor of $(m)_1$ in $(f^m)_3$.)  After substituting 
(2.9) into (2.2) for $f_j(t) = 1 + a_jt + b_jt^2$, $1 \le j \le 4$,
and taking common 
factors out of the rows, we find that
$$
W(m;t)= 
\gathered
\frac{m^3(m-1)}{12} t^6 
\vmatrix 1    &  \dots & 1 \\
 a_1   & \dots & a_4 \\
(m-1)a_1^2 + 2b_1 & \dots & (m-1)a_4^2 + 2b_4\\  
(m-2)a_1^3 + 6a_1b_1 &  \dots&(m-2)a_4^3 + 6a_4b_4 \endvmatrix.
\endgathered
\tag 2.10
$$
Let $W'(m)$ denote the determinant in (2.10); if $1 \neq m_0 \in
T(F)$, then  $W'(m_0) = 0$. 

\proclaim{Example 1 (Part Two)}
\endproclaim
We can dehomogenize (1.1) by setting $(x,y)\mapsto(1,t)$ and
normalizing so that $f_j(0) = 1$. This gives 
$$
\gathered
f_1(t) = 1 + \sqrt 2\ t - t^2,\qquad f_2(t) = 1 + i \sqrt 2\ t + t^2,\\
f_3(t) = 1 - \sqrt 2\ t - t^2,\qquad f_4(t) = 1 - i\sqrt 2\ t + t^2.
\endgathered
\tag 2.11
$$
Taking
$a_j = i^{j-1}\sqrt 2$ and $b_j = (-1)^j$ in (2.10), we obtain
$$
\gathered W'(m) =
\vmatrix 1   &1 &  1 & 1 \\
 \sqrt 2  & i\sqrt 2 & -\sqrt 2 & -i\sqrt 2 \\
2(m-2) &-2(m-2) & 2(m-2) & -2(m-2)\\  
2\sqrt2(m-5) &-2i\sqrt2(m-5) &-2\sqrt2(m-5)& 2i\sqrt2(m-5) \endvmatrix \\
= -128i(m-2)(m-5).
\endgathered
$$
Thus, $W'(m) = 0$ only for $m=2,5$. Remarkably, each root corresponds
to a linear dependence among the $\{f_j^m\}$. If we replace $\sqrt 2$
in (2.11) by a parameter $\mu>0$, then
the roots of $W'(m)$ work out to be  $m = 1
+ 2\mu^{-2}$, and $m = 2 + 6\mu^{-2}$. We shall show in Example 1
(Part Three) that, apart from $\mu = \sqrt 2$, these roots correspond to
linear dependence only when $(\mu,m) = (\sqrt 6, 3)$ and $(\sqrt{2/3},4)$. 

\head
3. $\ct(r,n,d)$ for small $r$, $n$,  $d$
\endhead
In this section, we describe
$\ct(r,n,d)$ in some simple cases.

\subhead
a. $(n,d)= (2,1)$
\endsubhead

It is very easy to show that  $\ct(2,1,r) =\{1,2, \dots, r-2\} $. If
$m \le r-2$, then $r > \binom {2+m\cdot 1-1}{2-1} = m+1$,
hence $m \in T(F)$ is forced. Suppose
$m = r-1$ and write $f_j(x,y) = \al_j x + \be_j y$ for $1 \le j \le r$.
Then the matrix giving $f_j^{r-1}$ with respect to the basis $\binom
{r-1}{\ell} x^{r-1-\ell}y^{\ell}$, for $0 \le \ell \le r-1$, has
Vandermonde determinant
$$
\vmatrix \al_1^{r-1}   & \al_1^{r-2}\be_1 & \dots  & \be_1^{r-1} \\
 \al_2^{r-1}   & \al_2^{r-2}\be_2 & \dots  & \be_2^{r-1} \\
 \vdots  & \vdots & \ddots & \vdots \\
\al_{r}^{r-1}   & \al_{r}^{r-2}\be_{r} & \dots  &
\be_{r}^{r-1} \endvmatrix = \prod_{i < j} \left( \al_i\be_j -
\al_j\be_i\right) \neq 0,
$$
since the $f_j$'s are pairwise non-proportional.
It follows that $r-1 \not \in T(F)$. If $m \ge r$, then take $f_j(x,y)
= \al_j x + \be_j y$ for $r+1 \le j \le m+1$ so that 
$F': = F \cup \{f_{r+1}, \dots, f_{m+1}\}$ consists of pairwise
non-proportional forms. Since $m \notin T(F')$,  $m \notin T(F)$.
\subhead
b. $(n,d)= (\ge 3,1)$
\endsubhead

It follows from the last section that
$\{1,\dots,r-2\} \in \ct(r,n,1)$ for $n \ge 3$. Are there
other tickets? We need a simple lemma.
\proclaim{Lemma 1}
Suppose 
$$
\sum_{j=1}^r \la_j\left(\al_{j1}x_1 + \dots + \al_{jn}x_n\right)^m =
0. \tag 3.1 
$$
Then for every $\ga = (\ga_1,\dots,\ga_n) \in \bold C^n$,
$$
\sum_{j=1}^r \la_j \left( \sum_{k=1}^n \al_{jk}\ga_k\right)
(\al_{j1}x_1 + \dots + \al_{jn}x_n)^{m-1} = 0. \tag 3.2 
$$
\endproclaim
\demo{Proof}
Apply the differential operator $\frac 1m\sum\limits_{k=1}^n \ga_k
\frac{\partial}{\partial x_k}$ to both sides of (3.1).
\qed
\enddemo
\proclaim{Lemma 2}

(i) If  $m \in T(F)$ for $F \in \cf(r,n,1)$, then $m-1 \in T(F)$.

(ii) If  $F \in \cf(r,n,1)$, then $\overline{m}(F) \le r-2$.

\endproclaim
\demo{Proof}
(i) If  $m \in T(F)$, then some equation (3.1) holds, where, say,
$\la_{j_0} \neq 0$. Since  $f_{j_0} \neq 0$,
there exists  $\ell$ so that $\la_{j_0}\al_{j_0\ell} \neq 0$. 
Now differentiate (3.1) with respect to $x_\ell$; (3.2) implies that
$m-1 \in T(F)$. 

(ii) Suppose (3.1) holds. We
claim by induction on $m$ that $m \le r-2$. If $m=1$, then the
assumption that the forms are pairwise non-proportional implies that
$r \ge 3$. Suppose now that the claim is true for $m-1$. There exists
$\ga$ so that $\sum_{k=1}^n \al_{rk}\ga_k = 0$, but $\sum_{k=1}^n
\al_{jk}\ga_k \neq 0$ (as $\al_j$ is not proportional to $\al_r$)  and
so the summand in (3.2) for $j=r$ disappears.
By induction, $m-1 \le r-3$, so $m \le r-2$.
\qed
\enddemo
If $r > \binom {n+m-1}{n-1}$, then $\{1,\dots,m\}$ is forced to be in
$T(F)$. Let $m(r,n)$ denote 
the largest integer $m$ so that $r > \binom {n+m-1}{n-1}$, and note
that $m(r,2) = r-2$. We have established the following theorem:
\proclaim{Theorem 2}
If $F \in \cf(r,n,1)$, then $T(F) = \{1,\dots,k\}$, where
$m(r,n) \le k \le m-2$.
\endproclaim

\proclaim{Example 3}
\endproclaim
Fix $(r,n)$ and let $m = m(r,n)$, so that $\binom {n+m}{n-1} \ge r >
\binom {n+m-1}{n-1}$. By Biermann's Theorem (see [11,\ p.31]),
$$
\{(i_1 x_1 + \dots + i_n x_n)^{m+1}: 0 \le i_k \in \bold Z,\ i_1 +
\dots + i_n = m+1\} \tag 3.3
$$
is a linearly independent set with $\binom{n+m}{n-1} \ge r$ elements. Let
$F$ be any subset of $r$ of the linear forms in (3.3). Then $m+1
\not\in T(F)$, but $\{1,\dots, m\}$ is forced to be in $T(F)$.
Thus $T(F)$ is precisely $\{1,\dots, m(r,n)\}$.
\smallskip
We strongly suspect that every ticket  allowed by Theorem 2 is
achieved, since the extremes appear.  We also strongly suspect that  a
proof of this result is not difficult.

An alternative presentation of the foregoing can be made using
``Serret's Theorem" (see [1] or [11,\ p.26]): given $\al_j \in \bold
C^n$, $1 \le j \le r$, $\al_j = (\al_{j1},\dots, \al_{jn})$, the
forms $\{ (\sum_k \al_{jk}x_k)^m\}$ are linearly independent if and
only if there exist  ``dual" forms $h_i(x_1,\dots,x_n)$ of degree $m$
so that $h_i(\al_j) = 0$ if $i \neq j$, but $h_i(\al_i) \neq
0$. Theorem 2 in [1] implies Lemma 3.2(ii).

\subhead
c. $r\le 3$
\endsubhead

These two cases are easily analyzed.
\proclaim{Theorem 3}

(i) $ \ct\{2,n,d\} = \{\emptyset\}$;

(ii) $\ct\{3,2,1\} = \{\{1\}\}$;

(iii) $ \ct\{3,n,1\} = \{\emptyset, \{1\}\}$, if $n \ge 3$;

(iv)  $ \ct\{3,n,d\} = \{\emptyset, \{1\}, \{2\}\}$, if $d \ge 2$.
\endproclaim
\demo{Proof}
First, if  $r=2$, $0 = \la_1f_1^m + \la_2f_2^m$ and
$\la \neq (0,0)$, then $f_1$ and $f_2$ are proportional.
If $r=3$ and $d=1$, then Theorem 2 and Example 3 imply (ii) and
(iii). In any event, Theorem 1 implies that $|T(F)| \le 1$,
``Fermat's Last Theorem" implies that $\overline{m}(F) \le 2$ and
(1.2)(iii) given $T(F) = \{2\}$.
\qed
\enddemo

\subhead
d. $(r,n,d) = (4,2,2)$
\endsubhead

The simplest remaining case is $(r,n,d) = (4,2,2)$. The heavy lifting
for this is done in the forthcoming  [12], 
which determines  all families $F \in \cf(4,2,2)$ for which
$3 \le m \in T(F)$. Since $4 >\binom{2+2-1}{2-1}$, 1 is forced, and 
$|T(F)| \le 3$ by Theorem 1. 

If $T(F)\subseteq \{1,2\}$, then the possibilities are $\{1\}$
or $\{1,2\}$. These are achieved, respectively by 
$$\{x^2,y^2,(x+y)^2,(x-y)^2\}, \qquad\{x^2,y^2,x^2+y^2,x^2-y^2\}.
$$

Suppose $m \in T(F)$ with $m \ge 3$.  Then  we may scale so that
$$
f_1^m + f_2^m = f_3^m +  f_4^m,\qquad m \ge 3.
\tag 3.4
$$
is an equation in binary quadratic forms.
Write the common sum in (3.4) as $p$, then $p \neq 0$,
$\{f_1^m,f_2^m\} \neq \{f_3^m,f_4^m\}$, and there do not exist 
$\al_j \in \bold C$ and $g \in \bold C[x,y]$ so that  $\al_1^m +
\al_2^m = \al_3^m + \al_4^m$ and $f_j = \al_jg$.
This is the equation studied in [12], where the uniqueness assertions
in the rest of this subsection  are established. (Uniqueness is taken
up to invertible linear changes of variable.)

\proclaim{Example 4}\endproclaim

For $m=3$, the solutions to (3.4) satisfying the non-triviality
conditions come from a single one-parameter family, first given, in a
different form, by J. Young in 1832 (see [3,\ p.554]).  
For $\al \notin \{0,1,-1\}$ and $\om = e^{2\pi i/3}$,  let
$$
\gathered
f_1(x,y) = \al x^2 - xy  + \al y^2, \qquad f_2(x,y) = -x^2 + \al xy -
y^2, \\
f_3(x,y) = \om \al x^2 -  xy  +  \om^2\al y^2, \qquad
f_4(x,y) = -\om x^2 + \al  xy -\om^2 y^2.
\endgathered
\tag 3.5
$$
Then 
$$
f_1^3(x,y) + \al f_2^3(x,y) = f_3^3(x,y) + \al f_4^3(x,y),
\tag 3.6
$$
and it is not hard to show that $T(\{f_1,f_2,f_3,f_4\}) = \{1,3\}$ for
every $\al$.
 (It is also true that $\om$ and $\om^2$ can be 
permuted in $f_3$ and $f_4$, giving a third pair of cubes in (3.6) with
the same sum.)
\smallskip
For $m=4$, there are two solutions to (3.4) satisfying the
non-triviality condition. 
\proclaim{Example 5}
\endproclaim
First, let
$$
\gathered
f_1(x,y) = x^2 + y^2,\qquad f_2(x,y) = \om x^2 + \om^2 y^2,\qquad 
f_3(x,y)= \om^2 x^2 + \om y^2,\\ \qquad
f_4(x,y) = xy.
\endgathered
\tag 3.7
$$
Then it is easy to verify that
$$
f_1 + f_2 + f_3 = 0, \quad f_1^2+f_2^2+f_3^2 = 6f_4^2, \quad
f_1^4+f_2^4+f_3^4 = 18f_4^4. 
$$
It follows from Theorem 1 that 
$T(\{f_1,f_2,f_3,f_4\}) = \{1,2,4\}$. The ``obvious" generalization of
(3.7) is presented in Example 8. 

After taking the linear transformation $(x,y) \mapsto (i(x-\om y),x -
\om^2y)$ in (3.7) and scaling, we obtain an integral 
version of $F$,
$$
\gathered
f_1(x,y) = x^2 + 2xy ,\qquad  f_2(x,y) = x^2 -  y^2, \\
f_3(x,y)= 2xy + y^2,\qquad f_4(x,y) = x^2+xy+y^2. 
\endgathered
\tag 3.8
$$

\smallskip

\proclaim{Example 6}
\endproclaim
The other solution to (3.4) for $m=4$ gives a family 
with ticket $\{1,4\}$:
$$
\gather
\sum_{\pm}(\sqrt 3\ x^2 \pm\sqrt 2\ xy - \sqrt3\ y^2)^4
= \sum_{\pm} (\sqrt 3\ x^2 \pm i \sqrt 2\ xy + \sqrt3\ y^2)^4; \tag
3.9
\endgather
$$
This example corresponds to $(\mu,m) = (\sqrt{2/3},4)$ in Example 1
(Part Two).
\smallskip
Finally, (1.1) gives the unique solution to (3.4) when $m=5$.

We summarize this discussion with a result whose proof relies on [12].
\proclaim{Theorem 4}
$$
\ct(4,2,2) = \{\{1\}, \{1,2\}, \{1,3\}, \{1,4\},
\{1,2,4\}, \{1,2,5\} \}.
$$
\endproclaim
As $d$ increases, so does $\ct(4,2,d)$.
Example 13 shows that $\{3\} \in \ct(4,2,4)$;
Example 14 shows that $\{4\} \in \ct(4,2,7)$.
\subhead
e. $\{(1,2,3)\} \notin \ct(4,n,d)$
\endsubhead

We close this section with a non-existence theorem. 
\proclaim{Lemma 3}
Suppose $\{f_j\} \in \cf(r,n,d)$ is two-dimensional; that is, there exist
$f$ and $g$, not proportional, so that $f_j = \al_j f + \be_j g$, $1
\le j \le r$. If $f_j' = \al_j x + \be_j y$, then $T(\{f_j\}) =
T(\{f_j'\})$. In particular,  $T(\{f_j\}) \subseteq \{1,\dots,r-2\}$.
\endproclaim
\demo{Proof}
One inclusion is obvious. For the other,
since no two $f_j$'s are proportional, the same can be said about the
 $(\al_j,\be_j)$'s. Suppose $\sum_j \la_jf_j^m = 0$ and 
let $P(x,y) =  \sum_j\la_j(\al_j x + \be_j y)^m$, so that $0 =
P(f,g)$. If $P= 0$, then $m \in T(\{f_{j}'\})$. Otherwise, $P$
splits into non-trivial linear factors, and $0 =  \prod_i(\ga_i f -
\de_i g)$ implies 
that $f$ and $g$ are proportional, a contradiction.
\qed
\enddemo

\proclaim{Theorem 5} 
For every $(n,d)$,  $ \{1,2,3\} \not\in \ct(4,n,d)$.
\endproclaim
\demo{Proof}
Let $F = \{f_1,f_2,f_3,f_4\}$ and  suppose there are three non-trivial
equations 
$$
\gather
0 = \la_{11}f_1 + \la_{12}f_2 + \la_{13}f_3 + \la_{14}f_4, \tag 3.10)(i \\
0 = \la_{21}f_1^2 + \la_{22}f_2^2 + \la_{23}f_3^2 + \la_{24}f_4^2, \tag
3.10)(ii \\ 
0 = \la_{31}f_1^3 + \la_{32}f_2^3 + \la_{33}f_3^3 + \la_{34}f_4^3. \tag
3.10)(iii
\endgather
$$ 
If the $f_j$'s were to satisfy a linear relation different from
(3.10)(i), then Lemma 3 would imply  $ 3 \not\in T(F)$. Also, Theorem
3 implies that $\la_{3j} \neq 0$ for all $j$.

Suppose first that one of $\la_{1k}$'s equals zero, say $\la_{14} =
0$. Then after scaling the $f_j$'s, we can assume that $f_1 + f_2 =
f_3$, and by 
Theorem 1 with $r=3$, $\la_{24} \neq 0$. Thus, after dividing by
$\la_{\ell4}$, and relabeling,
(3.10)(ii) and (3.10)(iii) imply that
$$
\gather
\be_{21}f_1^2 + \be_{22}f_2^2 + \be_{23}(f_1+f_2)^2 = f_4^2, \tag 3.11)(i \\
\be_{31}f_1^3 + \be_{32}f_2^3 + \be_{33}(f_1+f_2)^3 = f_4^3. \tag 3.11)(ii
\endgather
$$
If $\be_{21}u^2 + \be_{22}v^2 + \be_{23}(u+v)^2$ were a perfect
square, then (3.11)(i) would imply a linear relation among
$f_1,f_2$ and $f_4$, different from (3.10)(i), and as noted
above, this would be a contradiction. However, (3.11) implies that
$$
\left(\be_{21}f_1^2 + \be_{22}f_2^2 + \be_{23}(f_1+f_2)^2\right)^3 = 
\left(\be_{31}f_1^3 + \be_{32}f_2^3 +
\be_{33}(f_1+f_2)^3\right)^2; \tag 3.12
$$
since $f_1$ and $f_2$ are not proportional, (3.12) and the
argument of Lemma 3 imply that
$$
\left(\be_{21}u^2 + \be_{22}v^2 + \be_{23}(u+v)^2\right)^3 = 
\left(\be_{31}u^3 + \be_{32}v^3 +
\be_{33}(u+v)^3\right)^2. \tag 3.13
$$
But if $\Phi = F^3 = G^2$, then any irreducible factor of $\Phi$ must
occur to the sixth power. Thus, (3.13) implies that $\be_{21}u^2 +
\be_{22}v^2 + \be_{23}(u+v)^2$ is a square, a contradiction.

Therefore, we may assume that $\la_{1k} \neq 0$, and so after scaling, $f_4
= f_1 + f_2 + f_3$ and
$$
\gather
\be_{21}f_1^2 + \be_{22}f_2^2 + \be_{23}f_3^2 +
\be_{24}(f_1+f_2+f_3)^2 = 0 ,\tag 3.14)(i \\
\be_{31}f_1^3 + \be_{32}f_2^3 + \be_{33}f_3^3 +
\be_{24}(f_1+f_2+f_3)^3 = 0. \tag 3.14)(ii
\endgather
$$
By Bezout's Theorem, the intersection of the two curves
$$
\gather
G_2(u,v,w):= \be_{21}u^2 + \be_{22}v^2 + \be_{23}w^2 +
\be_{24}(u+v+w)^2 = 0 , \\
G_3(u,v,w):= \be_{31}u^3 + \be_{32}v^3 + \be_{33}w^3 +
\be_{34}(u+v+w)^3 = 0  
\endgather
$$
is at most six lines, unless the curves share a component. Since
$(f_1,f_2,f_3)$ 
parameterizes the intersection by (3.14), either $G_2\ | \ G_3$
or both share a  linear factor $\ell$. However, if $G_2 = \ell_1\ell_2$
for linear $\ell_j$, then $G_2(f_1,f_2,f_3) = 0$ implies
that $\ell_j(f_1,f_2,f_3) = 0$ for some $j$; this is different from $f_4 =
f_1+f_2+f_3$ and gives a contradiction as before.

The only remaining possibility is that  $G_2\ | \ G_3$ and $G_2$ is
irreducible. After relabeling once again, we see that this case is
$$
\gathered
\bigl(\si_1 u^2 + \si_2 v^2 + \si_3 w^2 +
2\si_4(uv+uw+vw)\bigr)\bigl(\ga_1 u + \ga_2 v + \ga_3 w\bigr) \\ 
= \ta_1 u^3 + \ta_2 v^3 + \ta_3 w^3 + \\
\ta_4(3u^2v+3uv^2+3uw^2+6uvw+3uw^2 + 3v^2w + 3vw^2).
\endgathered
\tag 3.15
$$
Since $\be_{3j} \neq 0$, we may assume in (3.15) that $\ta_4 \neq 0$
(hence  $(\ga_1,\ga_2,\ga_3) \neq (0,0,0)$) and that
$\ta_j \neq \ta_4$ for $j = 1,2,3$. 
By considering the coefficients of $u^2v, u^2w, uv^2,$ $ v^2w,
uw^2,vw^2$ and $uvw$ in (3.15), we find that
$$
\gathered
3\ta_4 = 2 \ga_1\si_4 + \ga_2\si_1 = 2 \ga_1\si_4 + \ga_3\si_1 =  
2 \ga_2\si_4 + \ga_1\si_2 = 2 \ga_2\si_4 + \ga_3\si_2 \\
= 2 \ga_3\si_4 + \ga_1\si_3 = 2 \ga_3\si_4 + \ga_2\si_3 =
(\ga_1+\ga_2+\ga_3)\si_4.
\endgathered
\tag 3.16
$$
Since $\ta_4 \neq 0$, it follows that $\si_4 \neq 0$. Further, 
it follows immediately from (3.16) that 
$$
0 = (\ga_2-\ga_3)\si_1 = (\ga_3-\ga_1)\si_2 = (\ga_1-\ga_2)\si_3. \tag 3.17
$$

If $\ga_1=\ga_2=\ga_3 := \ga \neq 0$, then (3.16)
reduces to 
$$
3\ta_4 = \ga(\si_1 + 2\si_4) = \ga(\si_2 + 2\si_4)= \ga(\si_3 +
2\si_4) = 3\ga\si_4.
$$
This implies that $\si = \si_2 = \si_3 = \si_4 := \si$, so $G_2(u,v,w)
= (u+v+w)^2$, a contradiction. 

If the $\ga_j$'s are not all equal, then after permutation of indices if
necessary, we may assume that $\ga_1 \neq \ga_2$ and $\ga_1 \neq
\ga_3$, hence $\si_2 = \si_3 = 0$ by (3.17). Then (3.16) becomes
$$
3\ta_4 = 2 \ga_1\si_4 + \ga_2\si_1 = 2 \ga_1\si_4 + \ga_3\si_1 =  
 2 \ga_2\si_4  = 2 \ga_3\si_4 = (\ga_1+\ga_2+\ga_3)\si_4.
$$
Thus $\ga_2 = \ga_3$, and so $2 \ga_2\si_4 = (\ga_1+\ga_2+\ga_3)\si_4$
implies  $\ga_1 = 0$. It then follows that $\si_1 = 2\si_4$, hence
$G_2(u,v,w) = 2\si_4(u+v)(u+w)$, a final contradiction. 
\qed
\enddemo
\head
4. A few interesting tickets 
\endhead

The examples in this section indicate some of the range of
$\ct$, and are inspired, for the most part, by the examples from [12].
The families mostly have $n=2$, consist of binomials or highly
symmetric trinomials, and rely on the properties of roots of
unity. 

Example 8 presents
highly dysfunctional families $F\in \cf(2v+2,2,2)$ with $|T(F)| = 3v$.
Example 9 gives families $F \in \cf(r,2,2)$ so that $T(F)$
contains $\{1,\dots,r-2\}$ as well as $r+s$, for any single integer $s$, $2
\le s \le r-1$. This is less dysfunctional than the previous example,
but with a larger jump. 
The families in Example 10 are not dysfunctional, but have an even
larger jump; for example, $\{1,2,8\} \in \ct(6,2,2)$. One special
case combines Examples 9 and 10 to show that  $\{1,2,3,4,8,14\} \in
\ct(6,2,2)$. Example 11 presents $F \in
\cf(a+2,2,a)$ so that $T(F)$ consists precisely of the divisors of
$a$. Example 12, essentially due to Molluzzo, shows that
$\overline{m}_r \ge \lfloor \frac{(r-1)^2}4\rfloor$ and also shows some
dysfunctionality. Finally, Examples 13  and 14 are variations
of two families due to L. Euler showing that $\{3\} \in \ct(4,2,4)$ and
$\{4\} \in \ct(4,2,7)$. 

We first need some elementary observations about roots of unity.
Write $\zeta_q = e^{\frac{2\pi i}q}$ for integral $q \ge 2$. 
We shall say that $\{f_0,\dots,f_{q-1}\}$ is {\it $q$-cyclotomic}
on $\{g_0,\dots,g_{q-1}\}$ if 
$$
f_{j} = \sum_{k=0}^{q-1} \zeta_q^{jk}g_k, \qquad 0 \le j \le q-1, \tag 4.1
$$
and the $g_k$'s involve disjoint sets of monomials. If $\{f_j\}$ is
$q$-cyclotomic on $\{g_k\}$, 
then $\{f_j^m\}$ is $q$-cyclotomic on $\{g_{m,k}\}$, where 
$$
g_{m,k} = \sum\Sb i_1 + \dots + i_m \equiv k \mod
q \\ 0 \le i_k \le q-1  \endSb g_{i_1}\dots g_{i_m}, \tag 4.2
$$
provided the $g_{m,k}$'s involve disjoint monomials. 
\proclaim{Lemma 4}
If $\{f_0,\dots,f_{q-1}\}$ is $q$-cyclotomic on
$\{g_0,\dots,g_{q-1}\}$, then 
$$
\text{Span}(f_0,\dots,f_{q-1}) = \text{Span}(g_0,\dots,g_{q-1}).
$$
In particular, $\{f_j\}$ is linearly dependent if and only if some $g_k = 0$.
\endproclaim
\demo{Proof}
By the orthogonality properties of roots of unity, (4.1) inverts to
$$
g_{\ell} = \frac 1q \sum_{j=0}^{q-1} \zeta_q^{-j\ell}f_j,  \qquad 0
\le \ell \le q-1.\tag 4.3 
$$
Since the $g_k$'s involve different monomials, they can only be
trivially dependent. 
\qed
\enddemo
\proclaim{Example 1 (Part Three)}
\endproclaim
Continuing the discussion in Example 1 (Part Two), suppose
$$
\gathered
f_1(t) = 1 + \mu t - t^2,\qquad f_2(t) = 1 + i \mu t + t^2,\\
f_3(t) = 1 - \mu t - t^2,\qquad f_4(t) = 1 - i \mu t + t^2.
\endgathered
\tag 4.4
$$
That is, $F$ is 4-cyclotomic on $\{1,\mu t,-t^2,0\}$.
To ensure pairwise non-proportionality in (4.4), we require $\mu \neq
0$. Suppose $1 < m \in T(F)$. Then it follows from Lemma 4 that
$g_{m,k} = 0$ for some $k$. But for $m > 1$, 
$$
\gathered
g_{m,0} = 1 + \dots,\qquad  g_{m,1} = m\mu t + \dots,\\
 g_{2v+1,2} = \dots +
t^{4v+2},\qquad g_{2v,3} = \dots - 2v\mu t^{4v-1},
\endgathered
$$
and since these $g_{m,k}$'s all have at least one non-zero term, the
only $g_{m,k}$'s which can vanish are $g_{2v,2}$ and $g_{2v+1,3}$. We have 
$$
\gathered
 g_{2,2}(t) = (\mu^2-2)t^2,\qquad g_{3,3}(t) = (\mu^3 - 6\mu)t^3,\\
 g_{4,2}(t) = (6\mu^2-4)(t^2 + t^6),\qquad
g_{5,3}(t) = (10\mu^3 - 20\mu)(t^3+t^7),\\ g_{6,2}(t) =
(15\mu^2-6)t^{2} + (\mu^6 - 30\mu^4+90\mu^2 - 20)t^6 +
(15\mu^2-6)t^{10}, \\
g_{7,3}(t) = (35\mu^3 - 42\mu)t^{3} + (\mu^7 -
42\mu^5+210\mu^3-140\mu)t^7 +  (35\mu^3 - 42\mu)t^{11}.
\endgathered
\tag 4.5
$$
We see from (4.5) that $g_{m,k} = 0$ for $m \le 7$ only if $(\mu,m) =
(\sqrt 2,2), (\sqrt 6, 3), (\sqrt{2/3},4)$ or $(\sqrt 2,5)$. (It is
evident that $m = 6,7$ yield nothing useful; Green's
Theorem implies that we need only check $m\le 8$; the computation of $g_{8,2}$
is left to the reader.)

We've seen that for  $\mu = \sqrt 2$, $F$ is equivalent to (1.1). 
(This includes the seemingly accidental $g_{2,2}\ | \ g_{5,3}$.)
If $\mu = \sqrt 6$, then after homogenizing and inverting by (4.2), we obtain
the identity
$$
\gathered 
6\sqrt 6(x^5y + xy^5) =  \\
(x^2 + \sqrt 6 xy - y^2)^3 - (-x^2  \sqrt 6 xy + y^2)^3
\\ = (ix^2 - \sqrt 6 xy +i y^2)^3 - (ix^2 + \sqrt 6 xy + iy^2)^3
\endgathered
\tag 4.6
$$
It happens that $6\sqrt 6(x^5y + xy^5)$ has four other representations
as a sum of two cubes, involving 24-th roots of unity; (4.6) is
derivable from (3.5) after an unpleasant change of variables and
choice of $\al$.
If $\mu = \sqrt{2/3}$, then after homogenizing and inverting, we
obtain  (3.9). 
\smallskip
\proclaim{Example 7}
\endproclaim
Let $f_j(x,y) = x+\zeta_q^jy$. Then $\{f_j^m\}$ is
$q$-cyclotomic 
on $\{g_{m,0},\dots,g_{m,q-1}\}$, where 
$$
g_{m,k}(x,y) = \sum \Sb i \equiv k \mod q \\ 0 \le i \le m \endSb
\binom mi x^{m-i}y^i.
$$
Thus  $g_{m,q-1}= 0$ if $m < q-1$, hence $\{f_j^m\}$ is
linearly dependent and so $m \in T(F)$, and $g_{m,k} \neq 0$ if $m \ge
q-1$, so $m \notin T(F)$.  This is just a special case of
Theorem 2.
\smallskip
\proclaim{Example 8}
\endproclaim
Suppose $q$ is odd and let  $f_j(x,y) = \zeta_q^j
x^2+\zeta_q^{-j}y^2$.  (Since $q$ is odd, no two $f_j$'s are proportional.)
Then  
$$
f_j^m(x,y) = \sum_{i=0}^m \binom mi \zeta_q^{(m-2i)j} x^{2m-2i}y^{2i},
$$
so 
$$
g_{m,k}(x,y) = \sum \Sb m-2i \equiv k \mod q \\ 0\le i \le m \endSb 
\binom  mi x^{2m-2i}y^{2i}. \tag 4.7
$$
Now let 
$$
F_q= \{f_0(x,y),\dots,f_{q-1}(x,y), xy\}.
$$
(Without $xy$, we are essentially in Example 7.)
By Lemma 4, 
$$
\text{Span}(f_0^m,\dots,f_{q-1}^m,(xy)^m) =
\text{Span}(g_{m,0},\dots,g_{m,q-1}, (xy)^m).
\tag 4.8
$$
It follows from
(4.8) that $m \in T(F_q)$ if and only if either $g_{m,i} = 0$ for some $i$ or
$g_{m,i}$ is a multiple of $(xy)^m$. 

In the first case, $m-2i \equiv k \mod q$ is
always soluble mod $q$, so if $m \ge q-1$, then there will exist $i$, $0
\le i \le m$, satisfying the congruence, and $g_{m,i} \neq 0$. On the
other hand, if $m \le q-2$, then  the 
equation $m-2i \equiv m+2 \mod q$ implies that $i \equiv q-1 \mod q$,
hence $g_{m,m+2} = 0$. Therefore, $T(F_q)$ contains
$\{1,\dots,q-2\}$. (Alternatively, $f_j(x,y) = g_j(x^2,y^2)$ for
linear $g_j$, and $m \le q-2$ is forced in $T(\{g_j\})$.)

In the second case, $(xy)^m$ can only occur in $g_{m,i}$ if $m$ is even,
$i = m/2$ and $k = 0$. Putting this into (4.7), we see that $(xy)^m$ is
the only term in $g_{m,0}$ if and only if $i = m/2$ is the only
solution to $2i \equiv m \mod q$ with $0 \le i \le m$, and this occurs
if and only if $m/2 < q$; that is, $m < 2q$. Thus $T(F_q)$ contains the
even integers $\{2, \dots, 2q-2\}$. Writing $q = 2v+1$, we see that
$r = |F_q| =  2v + 2$ and
$$
T(F_{2v+1}) = \{1,3,\dots,2v-1\} \cup \{2,4,\dots, 4v\} =
\{1,2,\dots,2v-1,2v,2v+2,\dots 4v\}.
$$
so that $|T(F_{2v+1})| = 3v = \frac 32 r - 3$. For example, (3.8) is
$F_3$. The explicit linear relations are not very interesting:
$$
\gathered
\sum_{j=0}^{q-1} \left( \zeta_q^j x^2+\zeta_q^{-j}y^2\right)^m =
0,\qquad \text{if $m < q$ is odd;} \\
\sum_{j=0}^{q-1} \left( \zeta_q^j x^2+\zeta_q^{-j}y^2\right)^m =
q\binom{m}{m/2}(xy)^m,\qquad \text{if $m < 2q$ is even.}
\endgathered
$$
\smallskip

\proclaim{Example 9}
\endproclaim
This next example is a variation on Example 8, giving a smaller
ticket, but one with a larger gap. We no longer assume that $q$ is
odd.  Let $\al \in \bold C$, and let 
$$
f_{j,q,\al}(x,y) = \zeta_q^jx^2+\al xy + \zeta_q^{-j}y^2, \qquad 0\le j \le
q-1, \tag 4.9
$$
and  let
$$
F_{q,\al} = \{f_{0,q,\al}(x,y), \dots , f_{q-1,q,\al}(x,y), xy\}.
$$
Since
$$
\left( \zeta_q^jx^2+\al xy + \zeta_q^{-j}y^2 \right)^m =
\sum\Sb a+b+c=m\\a,b,c\ge 0 \endSb \frac {m!}{a!b!c!}
\zeta_q^{j(a-c)}\al^bx^{2a+b}y^{b+2c},
$$
by setting $b = m - a -c$ above, we see that
$$
g_{m,k}(x,y) = \sum \Sb a-c \equiv k \mod q \\ 0\le a,c,m-(a+c) \endSb 
 \frac {m!}{a!c!(m-a-c)!}  \al^{m-a-c}x^{m+(a-c)}y^{m-(a-c)}. \tag 4.10
$$
Again, $m \in T(F_{q,\al})$ if and only if $g_{m,k} = 0$ or is a multiple of
$(xy)^m$. The complication here is that the coefficient of a monomial
in $g_{m,k}$ is a polynomial in $\al$, and may vanish when $\al$ is
suitably chosen.  If $m \le q-1$, then $a \equiv c \mod q$ and
$m \ge |a-c|$ imply that $a=c$ so that (4.10) becomes
$$
g_{m,0}(x,y) = \left(\sum_{a=0}^{\lfloor m/2 \rfloor}
 \frac {m!}{(a!)^2(m-2a)!}  \al^{m-2a}\right)x^my^m,
$$
and  $m \in T(F_{q,\al})$. Suppose now that $q \ge 3$ and $m = q + s$,
$2 \le s \le q-1$. We fix $k=0$. Then  $a-c \equiv 0 \mod q$ implies 
that $a - c \in \{-q,0,q\}$ and by (4.10),
$$
\gathered
g_{q+s,0}(x,y) = \\
\left(\sum_{a=0}^{\lfloor s/2 \rfloor}  \frac {(q+s)!}{a!(a+q)!(s -
2a)!}  \al^{s-2a}\right)(x^{s+2q}y^s + x^sy^{s+2q}) +  B(\al)(xy)^{s+q}. 
\endgathered
\tag 4.11
$$
Since $s\ge 2$, there exists $\al = \al_0$ which kills the coefficient
of $x^{s+2q}y^s + x^sy^{s+2q}$ in (4.11), so that $g_{q+s,0}$ is a
multiple of $(xy)^{s+q}$. In this case,
$$
 \{1,\dots,q-1,q+s\} \subseteq T(F_{q,\al_0}).
$$
We cannot rule out the possibility that  $T(F_{q,\al_0})$ has other
elements, if $\al_0$ is a root of other polynomials appearing in
different $g_{m.k}$'s. Unless that happens, 
we have $|F_{q,\al_0}| = q+1$ and if  $s = q-1$ is maximal, then
$\overline{m}({F}_{q,\al_0}) \ge q+s = 2q-1$. Thus,  
$\overline{m}_{r} \ge 2r-3$, and this is larger than
$\lfloor\frac{(r-1)^2}4 \rfloor$ for $r \le 8$. The stronger bound
$\overline{m}_6 \ge 14$ is shown in (4.16). 

The simplest concrete example is  $q=3$ and $m=5$. We have $\zeta_3 =
\om$ and (4.11) becomes
$$
g_{5,0} =  (5 + 10\al^2)(x^8y^2+x^2y^8).
$$ 
Taking  $\al_0 = \sqrt{-1/2}$, we define
$$
\gathered
f_{0,3,\al_0}(x,y) = x^2 + \al_0xy + y^2,\qquad f_{1,3,\al_0}(x,y) =\om x^2 +
\al_0xy + \om^2 y^2,  \\
f_{2,3,\al_0}(x,y) =\om^2 x^2 +\al_0xy + \om y^2, \qquad f_3(x,y) = xy. 
\endgathered
\tag 4.12
$$
Then $T(F_{3,\al_0})$ contains $\{1,2,5\}$.  After the linear change
$(x,y) \mapsto (i(y-cx),x+cy)$ for $c = \frac{\sqrt 6 - \sqrt 2}2 =
\sqrt{2-\sqrt 3}$, (4.12) becomes a scaled and permuted version of (1.1). 
\smallskip

\proclaim{Example 10}
\endproclaim
As a variation on the previous example, we again use (4.9), but
specify $q = 2v \ge 4$ to be even and let
$$
F_{q,\al}' = \{f_{0,q,\al}(x,y), \dots , f_{q-1,q,\al}(x,y)\} \in \cf(2v,2,2).
$$
Since $xy \notin F_{q,\al}'$, $m \in T(F_{q,\al}')$ only if
there exists $k$ so that $g_{m,k} = 0$. If  $m \le v-1$, then  $2v >
\binom{2m+1}1$ so $m$ is forced in $T(F_{q,\al}')$.  With a
careful choice of $\al$, another exponent can be added to the ticket. 
We are interested in the largest possible jump. If $m = 3v-1$, then by
(4.10),  
$$
\gathered
g_{3v-1,v}(x,y) = \\
\left(\sum_{a=0}^{v-1}  \frac {(3v-1)!}{a!(a+v)!(2v-1 - 2a)!}
\al^{2v-1-2a}\right)(x^{4v-1}y^{2v-1} + x^{2v-1}y^{4v-1}).  
\endgathered
\tag 4.13
$$  
Since $v \ge 2$, there  exists $\al = \al_0$ so that $g_{3v-1,v} =
0$; hence $\{1,...,v-1,3v-1\} \subseteq T(F_{2v,\al_0}')$. As
(4.13) gives an equation of degree $v-1$ satisfied by $\al^2$,  we
cannot expect ``nice" roots $\al_0$ for $v \ge 4$. By applying  the
inversion formula (4.3) to  $g_{3v-1,v} = 0$, we find that 
$$
0 = \sum_{j=0}^{2v-1} \zeta_{2v}^{-jv}f_{j,2v,\al_0}^m = 
\sum_{j=0}^{2v-1} (-1)^jf_{j,2v,\al_0}^m,
\tag 4.14
$$
hence the linear relation becomes two equal sums of $v$ $3v-1$-st
powers.

The simplest case is $v=2$, and (4.13) gives 
$$
g_{5,2}(x,y) = (10\al^3 + 20\al)(x^7y^3+x^3y^7),
$$
so $\al = \sqrt{-2}$. If we set $y \mapsto iy$, we get a scaled
version of (1.1). (The appearance of 2 in the ticket is a bonus.)
 If $v=3$, then 
$$
g_{8,3}(x,y) = 56\al(3 + 5\al^2 + \al^4)(x^{11}y^5+x^5y^{11}).
$$
By setting $\al_0 =  i \sqrt{\frac{5+ \sqrt{13}}2}$. 
and noting $\zeta_6 = - \om^2$, we see that (4.14) becomes, explicitly,
$$
\gathered
(x^2 + \al_0 xy + y^2)^8 + (\om x^2 + \al_0 xy + \om^2 y^2)^8 +
(\om^2 x^2 + \al_0 xy + \om y^2)^8 = \\
(x^2 - \al_0 xy + y^2)^8 + (\om x^2 - \al_0 xy + \om^2 y^2)^8 +
(\om^2 x^2 - \al_0 xy + \om y^2)^8 = \\
-3\bigl(\tfrac{1+\sqrt{13}}2\bigr)^4x^2y^2(4x^{12} - 13^{3/2}x^6y^6 +
4y^{12}).
\endgathered
\tag 4.15
$$
A  computer check of $g_{v,3v-1}$ for $v \le 40$ shows that it is
irreducible, except for $v=5$ which, amazingly, has a factor of $\al^2
+1$. Letting $\ep = \zeta_5$, and transposing as in (4.15), we make an
unexpected discovery:
$$
\sum_{j=0}^4(\ep^j x^2 + i x y + \ep^{-j}y^2)^{14}
= \sum_{j=0}^4(\ep^j x^2 - i x y + \ep^{-j}y^2)^{14}
= 5^7(xy)^{14} \tag 4.16
$$
We have inadvertently found another instance of Example 9, beyond the
range where it can be expected to work. The common sum
in (4.16), with $i$ replaced by $\al$, contains
 $(x^{19}y^9 + x^9y^{19})$ and $(x^{24}y^4 + x^4y^{24})$, each multiplied by
different  polynomials in $\al$ having $\al^2+1$ as a factor. 

If we now let $f_j = \ep^j x^2 + i x y + \ep^{-j}y^2$, $0 \le j \le 4$
and $f_5 = \sqrt{-5}\ xy$, and define  $F = \{f_0,\dots,f_5\}$, then 
$$
0 = \sum_{j=0}^4 f_j - \sqrt 5 f_5 = \sum_{j=0}^5 f_j^2 =
\sum_{j=0}^4 f_j^3 +\sqrt 5 f_5^3 =   \sum_{j=0}^5 f_j^4 =
\sum_{j=0}^5 f_j^8 = \sum_{j=0}^5 f_j^{14}.
\tag 4.17
$$
A computer check, using $\overline{m}_6 \le 24$, shows no other relations:
$T(F) = \{1,2,3,4,8,14\}$ and $\overline{m}_6 \ge 14$.

Finally, if we put $(x,y) = (1,-1)$ into (4.16) and use the closed form of
$\cos(\frac{2\pi k}5)$, we get a mysterious numerical identity:
$$
2\bigl( \sqrt 5 + 1 + 2i\bigr)^{14} + 2\bigl(-\sqrt 5 + 1 +
2i\bigr)^{14} + (4-2i)^{14} = 20^7.
$$
\smallskip
\proclaim{Example 11}
\endproclaim
This next example is embarrassingly simple by contrast. Fix $a$ and let
$$
\widehat F_a = \{x^a + y^a, x^a, x^{a-1}y, \dots , y^a\}.
$$
Then $\widehat F_a \in \cf(a+2, 2, a)$.  If $m \in T(\widehat F_a)$, then
there is a 
non-trivial equation
$$
 \sum_{k = 0}^a \la_k (x^{a-k}y^k)^m + \la_{a+1}(x^a + y^a)^m =0.
$$
But this is  impossible if $\la_{a+1} = 0$, so we may assume
 $\la_{a+1} = -1$ and
$$
\sum_{k = 0}^a \la_k (x^{a-k}y^k)^m = (x^a + y^a)^m = \sum_{\ell = 0}^m
\binom m\ell x^{ma - \ell a}y^{\ell a}. \tag 4.18
$$
But (4.18) occurs if and only if  $m \ | \ \ell a$ for each $\ell$, $0
\le \ell \le m$;
that is, if and only if $m$ is a divisor of $a$. Thus, $T(\widehat F_a)$
consists of the divisors of $a$ and $|T(\widehat F_a)|$
is the arithmetic function $d(a) \ll a$. This family is far from
dysfunctional, but by taking $a =
p^e$ for prime $p$,  we see that $T(\widehat F_{p^e}) =
\{1,p,p^2,\dots,p^e\}$. 
It follows that there is no bound on the ratio of
consecutive elements in $\ct$.
\smallskip
\proclaim{Example 12}
\endproclaim

We combine Examples 7 and 11. For parameters $a, q \ge 2$, let
$$
\widetilde F_{a,q} = \{x^a + y^a, x^a + \zeta_q y^a, \dots,  x^a +
\zeta_q^{q-1} 
y^a, x^a, x^{a-1}y, \dots, y^a\}.
$$
Then $\widetilde F_{a,q} \in \cf(q + a + 1,2,a)$. A dehomogenized version of
$\widetilde F_{a,q}$ is due to Molluzzo ([9], see [10,\ p.485]) in the
context that $qa \in T(\widetilde F_{a,q})$. By Lemma 4, $m \in
T(\widetilde F_{a,q})$ if and only if  there exists $k$, $0 \le k \le
q-1$, so that  
$$
\sum \Sb i \equiv k \mod q \\ 0\le i \le m \endSb 
\binom  mi x^{(m-i)a}y^{ia}
$$
is a linear combination of $\{x^{m(a-j)}y^{mj}: 0 \le j \le a\}$. 

For $k$, $0 \le k \le q-1$, let
$$
S_k(m;q,a) = \{ia:i \equiv k \mod q,\  0 \le i \le m\}
= \{(k + bq)a:  0 \le b \le \lfloor\tfrac{m-k}q\rfloor\}.
$$
and let $Z(m;q,a)$ denote the set of $k$ for which $S_k(m;a,q)$ only
contains multiples of $m$. 
Then  $m \in T(\widetilde F_{a,q})$ if and only if some $Z(m;a,q)$ is
non-empty. 

If $m \le q-1$, then $S_0(m;q,a) = \{0\}$, so  $m \in T(\widetilde  F_{a,q})$.
Since $S_0(m;q,a)$ consists of multiples of $qa$, if $m \ |\ qa$, then 
$0 \in Z(m;q,a)$, so  $m \in T(\widetilde F_{a,q})$. If $|S_k(m;q,a)|
\ge 2$ and  
$k \in  Z(m;q,a)$, then $m$ divides $ka$ and $(k+q)a$, so $m$ is a
divisor of $qa$. 

In the remaining case, $m \ge q$ is not a divisor of $aq$, but there
exists $k$ so that $m\ |\ ka$ and $S_k(m;q,a)$ is
a singleton, so that $k \ge m-(q-1)$. The existence of such $m$
depends on the arithmetic properties of $a$ and $q$. 
If $a=q=p$ is prime, then no such $m$ exists, and
$T(\widetilde F_{p,p}) = \{1,\dots,p,p^2\}$. On the other hand, if
$a=q=6$, say, then it is not hard to show that $m=8,10$ satisfy the
criteria and 
$$
T(\widetilde F_{6,6}) = \{1,2,3,4,5,6,8,9,10,12,18,36\}
$$
 Note that $F_{6,6}$ is (barely) dysfunctional. 

Since $a^2 \in T(\widetilde F_{a,a})$ where $ \widetilde F_{a,a} \in
\cf(2a+1,2,a)$, and $a(a+1) \in T(\widetilde F_{a+1,a})$ where $
\widetilde F_{a+1,a} \in \cf(2a+1,2,a)$, we obtain Molluzzo's bound:
$\overline{m}_r \ge \lfloor \frac{(r-1)^2}4 \rfloor$. 

\smallskip
\proclaim{Example 13}
\endproclaim
 The Euler-Binet (complete) solution to (3.5) for
$m=3$ (see [6,\ p.200] for a derivation) is
$$
\gather
f_1(x,y,z) = (x + 3y)(x^2 + 3y^2)z - z^4, \\
f_2(x,y,z) = (-x + 3y)(x^2 + 3y^2)z + z^4, \\
f_3(x,y,z) = (x^2 + 3y^2)^2 - (x - 3y)z^3, \\
f_4(x,y,z) =  -(x^2 + 3y^2)^2 + (x + 3y)z^3.
\endgather
$$
Here, $F = \{f_j\} \in \cf(4,3,4)$, and it can be shown that $T(F) =
\{3\}$. If we let $f_j'(x,y) = f_j(x,x,y)$, so $F' = \{f_j'\} \in
\cf(4,2,4)$, then $T(F') = \{3\}$ as well.
\smallskip
\proclaim{Example 14}
\endproclaim 
Euler's binary septic solution to (3.5) for $m=4$
(see [6,\ p.201] is 
$$
\gather
f_1(x,y) = x^7 + x^5y^2 - 2x^3y^4 + 3x^2y^5 + xy^6 ,\\
f_2(x,y) = x^6y - 3x^5y^2 - 2x^4y^3 + x^2y^5 + y^7,\\
f_3(x,y) = x^7 + x^5y^2 - 2x^3y^4 - 3x^2y^5 + xy^6 ,\\
f_4(x,y) = x^6y + 3x^5y^2 - 2x^4y^3 + x^2y^5 + y^7.
\endgather
$$
Here, $F = \{f_j\} \in \cf(4,2,7)$, and it can be shown that $T(F) =
\{4\}$. This can be seen by a shortcut, avoiding the machinery of
Theorem 1.  Suppose $\sum_{j=1}^4
\la_jf_j^m = 0$. Then consideration of the terms $x^{7m}$ and $y^{7m}$
shows that $\la_3 = -\la_1$ and $\la_4 = -\la_2$. Hence by
transposition, $\la_1(f_1^m - f_3^m) = \la_4(f_2^m -
f_4^m)$. Evaluation at (1,1) implies that $\la_1(4^m - (-2)^m) =
\la_4((-2)^m - 4^m)$, hence $\la_4 = -\la_1$, and the only possible
relation is $f_1^m + f_2^m - f_3^m - f_4^m = 0$. This equation can now
easily be checked for $m \le 8$.

It is apparently unknown whether there is a non-trivial solution to
the equation $f_1^4 + f_2^4 = f_3^4 + f_4^4$ in real binary forms of
degree $d$ for $3 \le d \le 6$. 
\smallskip
It seems clear that the vast, centuries-old literature of Diophantine
parameterizations can be data-mined for many more interesting examples
of tickets. 

\head
5. Generalizations, speculations and   open questions 
\endhead

This final section collects a number of miscellaneous remarks.

There is no reason to restrict our attention to forms over $\bold C$;
for any field $k$, one might just as easily define $\cf_k(r,n,d)$,
where the forms have coefficients in $k$. It would be particularly
interesting to study these for $k = \bold Q$ or $\bold R$. It is
proved in [12] that there is no 
solution to (3.4) for real polynomials when $m=5$. Thus, $(1,2,5)
\notin \ct_{\bold R}(4,2,2)$, but (3.8) shows that $(1,2,4)
\in \ct_{\bold Q}(4,2,2)$. The forms in the highly dysfunctional
family $F_q$ in Example 8 are $\bold R$-combinations of $x^2+y^2$,
$i(x^2-y^2)$ and $xy$, and hence can be made real by taking $(x,y)
\mapsto (x+iy,x-iy)$. It seems hard to decide
whether there exists a change of variables which makes a given family
real, or rational.

One can also sharpen the definition of $T(F)$. Let 
$$
\de_m(F) = |F| - \text{dim(span(}\{f_j^m\}),
$$
so that $T(F) = \{m: \de_m(F) > 0\}$. 
\proclaim{Conjecture 2}
$$
\sum_{m=1}^\infty \de_m(F) \le \binom {r-1}2. \tag 5.1
$$
\endproclaim
This conjecture implies Theorem 1 of course, and is valid for the
examples given where 
$|T(F)| = \binom{r-1}2$ for $r = 3,4$. Further, if $(n,d) = (2,1)$, 
then $\de_m(F) = r - (m+1)$ for $m \le r-2$, and  (5.1) is sharp. 
 
As my algebraic geometry friends have pointed out, the topic of this
paper is ``really" the study of
rational curves lying on the intersection of several
Fermat varieties of different degrees. The rub in proving Conjecture 1
would be trying to ensure that the curve doesn't lie on other Fermat
varieties. By Theorem 2, every finite set of positive integers $A$ is
a subset of a member of $\ct(r,2,1)$ for $r$ sufficiently large. What
makes a ticket interesting is as much what's {\it not} in it as well as
what's in it. We can weaken Conjecture 1 enough to make it very
plausible: suppose $A = \{m_k\}$ is finite and $w \notin A$. Does there exist 
$F$ so that  $A \subseteq T(F)$ and $w \notin T(F)$?
If not, then the linear dependence of  $\{f_j^{m_k}\}$  forces the
dependence of  $\{f_j^w\}$. This seems unlikely.

\proclaim{Example 1 (Part Four)}
\endproclaim

 Consider the intersection 
of the plane $M_1(t) = 0$ and the conic  $M_2(t)
= 0$. Then $t_4 = -(t_1 + t_2 + t_3)$, and 
$$
t_1^2 + t_2^2 + t_3^2 + (-(t_1 + t_2 + t_3))^2 = 
2(t_1^2 + t_2^2 + t_3^2 + t_1t_2 + t_1t_3 + t_2t_3) = 0. \tag 5.2
$$
It is not difficult to diagonalize (5.2)  as
$$
(t_1 - t_2)^2 + 2(t_1 + t_2)^2 +  (t_1 + t_2 + 2t_3)^2 = 0. \tag 5.3
$$
Let
$$
t_1 - t_2 = x^2 - y^2, \quad\sqrt 2(t_1 + t_2) = 2xy, \quad 
i(t_1 + t_2 + 2t_3) = x^2 + y^2
$$ 
be the usual Pythagorean parameterization of (5.3). After solving for
the $t_i$'s, we discover that $t_j = f_j$ (c.f. (1.1))!

The fact that $M_5(t)= 0$ as well can be explained in
several ways. First, 
$$
\gathered
t_1^5 + t_2^5 + t_3^5 + (-(t_1 + t_2 + t_3))^5 = \\
-5(t_1+t_2)(t_1+t_3)(t_2 + t_3)(t_1^2 + t_2^2 + t_3^2 + t_1t_2 +
t_1t_3 + t_2t_3), 
\endgathered
\tag 5.4
$$
so that (5.2) implies that (5.4) vanishes as well. 

More generally, let $e_k$, $1 \le k \le 4$, denote the
$k$-th elementary symmetric function in four variables.
By Newton's Theorem, if $p(t_1,t_2,t_3,t_4)$ is a symmetric
polynomial, then there is a polynomial $P$ so that $p =
P(e_1,e_2,e_3,e_4)$. If $p$ is {\it any} symmetric form of degree
five in four variables (not just $M_5$), then 
$$
p = \al_1 e_1^5 + \al_2 e_1^3e_2 + \al_3 e_1^2e_3 + \al_4 e_1e_2^2 +
\al_5 e_1e_4 + \al_6e_2e_3. \tag 5.5
$$
If $M_1(t)=M_2(t) = 0$, then $e_1(t) = e_2(t) = 0$, and so
$p = 0$ as well; an explicit representation of $p \in (e_1,e_2)$ is
found by solving for $\al_k$ in (5.5).

This construction generalizes. Suppose $r \ge 4$ and let $A_r$ denote
the set of integers which {\it cannot} be written in the form $a(r-1)
+ br$ for non-negative integers $a$ and $b$. It is well-known that
$|A_r| = \binom 
{r-1}2$ and the largest element in $A_r$ is $r(r-1)-r-(r-1)$. Suppose
$$
t_1^k + \dots + t_r^k = 0, \qquad 1 \le k \le r-2, \tag 5.6
$$
then the same argument given above implies that any homogeneous
symmetric polynomial  in $r$ variables whose degree lies in
$A_r$ will vanish. If we could find pairwise non-proportional
polynomials $f_j$, $1 \le j \le  r$, satisfying
(5.6), then we would thereby construct a family $F$ of $r$ polynomials
with $T(F) = A_r$.  However, it appears
that no such parameterization of (5.6) exists for $r \ge 5$.

There is another,  geometric, interpretation of (1.1).  Define the four
pairs of complex numbers  $(\be_j,\ga_j)$ by the factorization
$f_j(x,y) = \sigma_j(x - \be_jy)(x-\ga_j y)$. Under the standard
stereographic projection of $\bold C$ to the unit sphere $S^2$, the
four pairs $(\be_j,\ga_j)$ map to the antipodal pairs of vertices of
an inscribed cube. This is reminiscent of Klein's classical work on
the icosahedron.
\smallskip
We conclude this section with some more open problems.

Example 9 shows that $|T(F)|$ can be as large as $\frac 32 r -
3$. This is probably not maximal. Is $|T(F)| = {\Cal O}(r)$? If not,
what is the true growth rate? Can $c\cdot r^2$ be achieved for some $c
> 0$? Does $\lim\limits_{r\to\infty} r^{-2}\overline{m}_r$ exist?

Is it true that $T(F)$ contains at most $r-2$ consecutive
integers? This is true in all the dysfunctional families presented
here, and would subsume Theorem 5.

Finally, in view of the last three remarks in $\S1$, we make the
following definitions:
$$
\gather
\ct(r,n,\infty) = \bigcup_d \ct(r,n,d) = \lim_{d\to\infty} \ct(r,n,d); \\
\ct(r,\infty,d) = \bigcup_n \ct(r,n,d) = \lim_{n\to\infty} \ct(r,n,d); \\
\ct(r,\infty,\infty) = \bigcup_{n,d} \ct(r,n,d).
\endgather
$$
Since  $\ct(r,\infty,\infty) $ is finite, it must actually be achieved
by a finite  $\ct(r,n,d)$. Can we compute 
 ``minimal" $d(n,r)$ so that $\ct(r,n,d(n,r)) =
\ct(r,n,\infty)$, ``minimal" $n(d,r)$ so that $\ct(r,n(d,r),d) =
\ct(r,\infty,d)$ or ``minimal" $(n(r), d(r))$ so that $\ct(r,n(r),d(r)) =
\ct(r,\infty,\infty)$? 

\head
6. Acknowledgments
\endhead

I thank my colleagues Mike Bennett, Sean Sather-Wagstaff
and Jack Wetzel, the 
organizers of the Analytic Number Theory, Commutative Ring Theory and
Geometric Potpourri seminars at UIUC, for their tolerance in letting me
woodshed this material in their groups. I  also thank Ricky
Pollock, Marie-Francoise Roy and Micha Sharir for inviting me to speak at the
DIMACS workshop on Algorithmic and Quantitative Aspects of Real
Algebraic Geometry in Mathematics and Computer Science in March 2001.

I thank my colleagues Nigel Boston,  Dan Grayson  and Marcin Mazur for useful
conversations. Marcin alerted me to the work of Mark
Green. I would also like to thank J\'anos Koll\'ar, Noam Elkies, Mark
Green and 
Gary Gundersen  for helpful email correspondence.

More than is customary in such acknowledgments, the author wants to emphasize
that he is solely responsible for any errors in content or
aberrations in mathematical taste. 
 
\enddocument
\end